\documentclass[11pt]{article}

\usepackage{amsmath}
\usepackage{amsfonts}
\usepackage{amssymb}
\usepackage{amsthm}
\usepackage{graphics}
\usepackage{amscd}
\usepackage{graphicx,epsfig}
\usepackage{color}
\usepackage[all]{xy}
\usepackage{url}

\DeclareMathOperator{\Div}{div}
\DeclareMathOperator{\im}{im}

\DeclareMathOperator{\argth}{argth}

\renewcommand{\epsilon}{\varepsilon}

\newcommand{\boH}{\mathcal{H}}

\newcommand{\R}{\mathbb{R}}
\newcommand{\Q}{\mathbb{Q}}

\newcommand{\C}{\mathbb{C}}

\renewcommand{\H}{\mathbb{H}}

\newcommand{\N}{\mathbb{N}}
\newcommand{\E}{\mathbb{E}}

\newcommand{\eps}{\varepsilon}
\newcommand{\der}[2]{\dfrac{\partial #1}{\partial #2}}
\newcommand{\dd}{\mathrm{d}}

\newcommand{\Ome}{\Omega}

\newtheorem{defn}{Definition}
\newtheorem{thm}{Theorem}

\newtheorem{prop}[thm]{Proposition}
\newtheorem{lem}[thm]{Lemma}

\newcommand{\barre}[1]{\overline{#1}}
\renewcommand{\phi}{\varphi}
\newcommand{\dis}{\displaystyle}

\newtheorem*{thm*}{Theorem}

\newtheorem*{claim*}{Claim}

\theoremstyle{remark}

\newtheorem{remarq}{Remark}
\newtheorem*{rem*}{Remark}

\newcounter{remark}

\newcounter{case}

\newcounter{construction}

\newcounter{fact}

\newcounter{step}

\newcommand{\wtilde}{\widetilde}
\newcommand{\Nil}{\mathrm{Nil}_3}

\title{The half space property for cmc $1/2$ graphs in $\E(-1,\tau)$}
\author{Laurent Mazet\thanks{The author was partially supported by the
ANR-11-IS01-0002 grant.}}
\date{}

\begin{document}

\maketitle

\begin{abstract}
In this paper, we prove a half-space theorem with respect to constant mean
curvature $1/2$ entire graphs in $\E(-1,\tau)$. If $\Sigma$ is such an entire
graph and $\Sigma'$ is a properly immersed constant mean curvature $1/2$ surface
included in the mean convex side of $\Sigma$ then $\Sigma'$ is a vertical
translate of $\Sigma$. We also have an equivalent statement for the non mean
convex side of $\Sigma$.
\end{abstract}

\section{Introduction}
In the theory of constant mean curvature surfaces, the half-space property is a
problem that have received many contributions in recent years. If $\Sigma$ is a
properly embedded constant mean curvature $H_0$ surface in a Riemannian
$3$-manifold $M$, we may wonder about the existence of an other properly
embedded constant mean curvature $H_0$ surface $\Sigma'$ which has no
intersection with $\Sigma$. If such a surface exists, can we say anything about
its geometry?

One of the first results about this question is the half-space theorem of Hoffman
and Meeks \cite{HoMe}. It says that a minimal surface of $\R^3$ on one side of a
plane $P$ is a plane parallel to $P$. The same type of results have been proved
by several authors in other homogeneous spaces (see \cite{RodRos,RoRo,HaRoSp,
DaHa,DaMeRo,Pen,Maz15,RoScSp}). In most of these results, the half-space
property is studied with respect to a surface $\Sigma$ which is parabolic. The
parabolicity is really an important hypothesis as it has been proved in
\cite{Maz15} and \cite{RoScSp}.

The first result involving non parabolic surfaces is due to Daniel, Meeks and
Rosenberg \cite{DaMeRo} and concerns minimal surfaces in the Heisenberg space
$\Nil$. As a Riemannian homogeneous space, $\Nil$ is a Killing Riemannian
submersion over $\R^2$. So, in $\Nil$, we can consider surfaces called graphs
that are images of sections of the submersion. Among these surfaces, the minimal
entire graphs (graph over the whole $\R^2$) have been classified by Fernandez
and Mira \cite{FeMi}; certain ones are parabolic and others are not. Actually,
Daniel, Meeks and Rosenberg proved that if $\Sigma'$ is a properly immersed
minimal surface in $\Nil$ on one side of an entire minimal graph $\Sigma$ then
$\Sigma'$ is a vertical translate of $\Sigma$.

An other situation where interesting non parabolic surfaces appears is in
considering entire constant mean curvature $1/2$ graphs in $\H^2\times\R$. In
fact, Daniel and Hauswirth \cite{DaHa} proved that there is an isometric
correspondence between such surfaces and entire minimal graphs in $\Nil$. This
suggests that a half-space property similar to the one of $\Nil$ should be true
for these surfaces in $\H^2\times\R$. Actually, we have partial results in this
direction, the half-space property had been established by Nelli and Sa Earp
\cite{NeSa} for a particular entire graph which is rotationally symmetric. More
recently Cartier and Hauswirth \cite{CaHa} proved also the half-space property
for a family of entire graphs that can be obtained by deforming the rotationally
symmetric one.

In this paper, we prove the half-space property for all entire cmc $1/2$ graphs
in $\H^2\times\R$. In fact, we prove this for entire cmc $1/2$ graphs of
$\E(-1,\tau)$. $\E(-1,\tau)$ denotes the family of simply connected homogeneous
spaces which are Killing Riemannian submersions over $\H^2$; we have
$\H^2\times\R=\E(-1,0)$. Our main theorem (Theorem~\ref{th:half}) is then
\begin{quote}
\textbf{Theorem.} Let $\Sigma$ be an entire constant mean curvature $1/2$ graph
in $\E(-1,\tau)$ and $\Sigma'$ be a properly immersed constant mean curvature
$1/2$ surface in $\E(-1,\tau)$. If $\Sigma'$ is included in the mean convex side
of $\Sigma$ then $\Sigma'$ is a vertical translate of $\Sigma$.
\end{quote}

We also have a statement when $\Sigma'$ is included in the non-mean convex side
of $\Sigma$.

Our strategy of proof is similar to the one of \cite{DaMeRo} and \cite{RoScSp}
which consists in constructing a family of barriers that converges to a vertical
translate of $\Sigma$. The main difficulty is to prove that the barriers we
construct actually converge to a vertical translate of $\Sigma$. So the main
part of the proof (Proposition~\ref{prop:uniq}) is devoted to a
uniqueness result for the exterior Dirichlet problem associate to the constant
mean curvature $1/2$ equation. If $u_0\le u_1$ are two solutions of our exterior
Dirichlet problem, we construct $(u_t)_{0\le t\le 1}$ a family of solutions of
the exterior problem that goes continuously from $u_0$ to $u_1$. Then
$\partial_t u_t$ defines a Jacobi field on the graph of $u_t$. We use the
associated family introduced by Daniel in \cite{Dan2} to study this Jacobi field
on a minimal surface in $\Nil$. In fact, we prove uniqueness of this Jacobi
field. In some sense, it can be interpreted as an infinitesimal version of the
half-space theorem of Daniel, Meeks and Rosenberg. Using this uniqueness, we can
reintegrate $\partial_t u_t$ in order to prove that $u_0$ and $u_1$ should
differ by a constant which leads to the conclusion.

In Section~\ref{sec:def}, we recall some definitions about the $\E(\kappa,\tau)$
spaces and explain what are the entire constant mean curvature graphs in them.
Section~\ref{sec:estim} is devoted to two gradient estimates for solutions of
the constant mean curvature equation, these estimates are used in the
construction of the barriers and in the proof of the uniqueness result. In
Section~\ref{sec:half}, we prove our main theorem assuming the uniqueness result
(Proposition~\ref{prop:uniq}). This uniqueness result is proved in
Section~\ref{sec:unik}. Finally, in Appendix~\ref{app:a}, we give some
computations concerning Killing Riemannian submersions and, in
Appendix~\ref{app:b}, we construct two barriers that are used in
Section~\ref{sec:half}.

\section{Entire cmc $1/2$ graphs in $\E(-1,\tau)$}
\label{sec:def}
\subsection{The ambient spaces $\E(\kappa,\tau)$}

In this section, we give a quick introduction to the ambient space
$\E(\kappa,\tau)$ when $\kappa\le 0$; for a complete description we refer to
\cite{Dan2}. 

The space $\E(\kappa,\tau)$ is a simply connected homogeneous space with an
isometry group of dimension at least $4$. For $\kappa\le 0$, we define
$D_\kappa=\{(x,y)\in\R^2\,|\,1+\kappa(x^2+y^2)\ge 0\}$ ($D_0=\R^2$). A
model for the ambient space $\E(\kappa,\tau)$ is then $D_\kappa\times\R$ with
the following complete Riemannian metric
$$
\dd s_{\kappa,\tau}^2=\lambda_\kappa^2(\dd x^2+\dd y^2)+(2\tau\lambda_\kappa
(y\dd x-x\dd y)+\dd z)^2
$$
where
$$
\lambda_\kappa=\frac{2}{1+\kappa(x^2+y^2)}
$$

The map $\pi:\E(\kappa,\tau)\rightarrow (D_\kappa,\lambda_\kappa^2(\dd x^2+\dd
y^2)), (x,y,z)\mapsto(x,y)$ is then a Killing
Riemannian submersion (see definitions in Appendix~\ref{app:a}) over either the
hyperbolic space $\H^2(\kappa)$ of curvature $\kappa$ or the Euclidean space
$\R^2$ ($\kappa=0$). Moreover the unit Killing vector field is $\xi=\partial_z$.

If $\mu>0$, the map $h_\mu: (x,y,z)\mapsto \mu(x,y,z)$ is a diffeomorphism from
$\E(\kappa,\tau)$ to $\E(\frac\kappa{\mu^2}, \frac\tau\mu)$ such that
$$
h_\mu^*(\dd s_{\frac\kappa{\mu^2},\frac\tau\mu}^2)=\mu^2\dd s_{\kappa,\tau}^2
$$
So the study of the geometry of the spaces $\E(\kappa,\tau)$ when $\kappa<0$
reduces to the one of the spaces $\E(-1,\tau)$; so in the following we only
focus to these spaces. We will denote $\H^2=\H^2(-1)$ and
$\lambda=\lambda_{-1}$, we also use $\barre \nabla$ to denote the Levi-Civita
connection on $\E(-1,\tau)$.

The vector field $\xi$ generates a flow $(\phi_t)_t$, in the following we will
denote by $p+t$ the point $\phi_t(p)$ where $p\in\E(-1,\tau)$. This notation is
coherent with the identification of sections with functions in the model.

The coordinate $z$ defines a function on $\E(-1,\tau)$, in the following we will
consider the gradient of this function. So we introduce the vector field $\zeta$
with the following expression
\begin{equation}\label{eq:zeta}
\zeta=\overline\nabla z=-2\tau y F_1+2\tau xF_2+\xi
\end{equation}
where $F_1=\frac1\lambda \partial_x-2\tau y\partial_z$ and $F_2=\frac1\lambda
\partial_y+2\tau x\partial_z$ are the horizontal lift of an orthonormal frame of
$\H^2$ (let us notice that $(F_1,F_2,\xi)$ is an orthonormal frame of
$\E(-1,\tau)$).


\subsection{The mean curvature of graphs in $\E(-1,\tau)$}

In $\E(-1,\tau)$, a surface is called a graph if it is the image of a smooth
section $\sigma : \Ome\subset\H^2\rightarrow \E(-1,\tau)$, this image is called
the graph of $\sigma$ (see Appendix~\ref{app:a}). For such a section, we can
define a vector field $G\sigma$ on $\Ome$ by the following property:
$$
(G\sigma(p),X)_{\H^2}=(\xi,\dd_p \sigma(X))_{\E(-1,\tau)}\textrm{ for any }X\in
T_p\H^2
$$
When $\tau=0$ and $\sigma$ is a function, $G\sigma$ is the gradient of
$\sigma$; in general, $G\sigma$ will play the role of the gradient.

In fact, the graph of $\sigma$ has constant mean curvature $H_0$ if
\begin{equation}\label{eq:cmc}
\Div_{\H^2}\left(\frac{G\sigma}{\sqrt{1+\|G\sigma\|^2}}\right)=2H_0
\end{equation}
where the mean curvature is computed with respect to the upward pointing normal
(see \eqref{eq:appcmc} in Appendix~\ref{app:a} for a proof of \eqref{eq:cmc}).

If we use the coordinates given by the model, a section $\sigma$ can be
identified to a function (\textit{i.e.} $z=\sigma(x,y)$) and $G\sigma$ has the
following expression:
$$
G\sigma=\frac1{\lambda^2}(\sigma_x+2\tau\lambda y)\partial_x+
\frac1{\lambda^2}(\sigma_y-2\tau\lambda x)\partial_y
$$


\subsection{Entire cmc $1/2$ graphs in $\E(-1,\tau)$}
\label{sec:entiregraph}

If $\sigma$ is a section defined on the whole $\H^2$ and its mean curvature is
constant $H_0$, it is known that this mean curvature is less than $1/2$ (see
\cite{Pen2}). In fact, we are interested in entire graph with mean curvature
$1/2$; so we consider sections $\sigma$ defined on the whole $\H^2$ which are
solution of 
\begin{equation}\label{eq:cmc12}
\Div_{\H^2}\left(\frac{G\sigma}{\sqrt{1+\|G\sigma\|^2}}\right)=1
\end{equation}
Such an entire graph $\Sigma$ bounds two connected components of $\E(-1,\tau)$,
one is above the graph, it is the mean convex side of $\Sigma$, and one is
below.

In fact, the space of all constant mean curvature $1/2$ entire graphs in
$\E(-1,\tau)$ is classified (see \cite{DaHa}, \cite{HaRoSp} and \cite{Pen}). Let
us give some explanation on the classification.

Let $X:\Sigma\rightarrow \E(-1,\tau)$ be a simply connected constant mean
curvature $1/2$ surface in
$\E(-1,\tau)$. On this surface, we can consider four geometric data: they are its
metric $g$, its shape operator $S$, a function $\nu$ and a vector field $T$
such that along $\Sigma$ the vector field $\xi$ can be written $\xi=T+\nu N$
where $N$ is the unit normal vector to $\Sigma$. In \cite{Dan2}, Daniel prove
that these data satisfies to equations which are necessary and sufficient for
the existence of the immersion $X$.

Now if we consider $\theta$ such that $\tau+\frac i2= e^{i\theta}
\sqrt{\tau^2+\frac14}$. We can consider the following data
\begin{align*}
g'&=g\\
S'&=e^{\theta J}(S-\frac 12I)\\
\nu'&=\nu\\
T'&=e^{\theta J}T
\end{align*}
These new data solve the conditions introduced by Daniel \cite{Dan2} for the
existence of a minimal immersion $X'$
of $\Sigma$ in $\E(0,\tau')$ with $\tau'=\sqrt{\tau^2+\frac14}$. Actually this
gives an isometric correspondence between minimal immersions in $\E(0,\tau')$ and
cmc $1/2$ immersions in $\E(-1,\tau)$. By the work of Daniel, Hauswirth,
Rosenberg, Spruck and Penafiel \cite{DaHa, HaRoSp, Pen}, we know that being an
entire graph is preserved by this correspondence. So every entire cmc $1/2$
graph in $\E(-1,\tau)$ corresponds to an entire minimal graph in $\E(0,\tau')$.
By the work of Fernandez and Mira \cite{FeMi}, we know that the space of entire
minimal graphs in $\E(0,\tau')$ is in correspondence with the space of quadratic
holomorphic differential over $\C$ or the unit disk. This gives the
classification of entire cmc $1/2$ graphs in $\E(-1,\tau)$

\section{Gradient estimates}
\label{sec:estim}
In the sequel, we need several gradient estimates for constant mean curvature
$1/2$ graphs in $\E(-1,\tau)$. Actually, we consider these surfaces as graphs of
functions in the model; in consequence the estimates we get depend on our choice
of the model.


\subsection{Some preliminary computations}

Let $\sigma$ be a section of $\E(-1,\tau)$ whose graph has constant mean
curvature $H_0$, $\sigma$ is a solution of \eqref{eq:cmc}. Let $N$ denote the
upward pointing unit
normal to the graph $\Sigma$ of $\sigma$. On $\Sigma$, we consider the function
$\nu=(N,\xi)$. Since $\xi$ is a unit killing vector field, $\nu$ is a Jacobi
function on $\Sigma$ so: 
$$
\Delta_\Sigma \nu =-(Ric(N,N)+\|S\|^2)\nu.
$$
where $Ric$ is the Ricci tensor of $\E(-1,\tau)$. Actually, we have
$Ric(N,N)=-(1+2\tau^2)+\nu^2(1+4\tau^2)$ (see \cite{Dan2}), thus we get:
$$
\Delta_\Sigma \nu =-(-(1+2\tau^2)+\nu^2(1+4\tau^2)+\|S\|^2)\nu.
$$

We also have
$$
\Delta_\Sigma \frac1\nu=-\frac{\Delta_\Sigma\nu}{\nu^2}+
2\frac{\|\nabla_\Sigma\nu^2\|^2}{\nu^3}.
$$
So if we introduce the operator $Lu=\Delta_\Sigma
u-2\nu(\nabla_\Sigma\frac1\nu,\nabla_\Sigma u)$, we have:
$$
L\frac1\nu=(-(1+2\tau^2)+\nu^2(1+4\tau^2)+\|S\|^2)\frac1\nu\ge
-\frac{1+2\tau^2}\nu.
$$
Using the model, we denote by $h$ the restriction of the $z$ coordinate to
$\Sigma$. We then have
$$
\nabla_\Sigma h= \zeta^\top
$$
where $X^\top$ denotes the orthogonal projection of $X$ on $T\Sigma$.
Using the expression~\eqref{eq:zeta}, we obtain the following estimate:
\begin{equation}\label{eq:est1}
\|\nabla_\Sigma h\|=\|\zeta\|^2-(\zeta,N)^2\ge 1-\nu^2-c_1\nu
\end{equation}
where $c_1$ is a positive constant. For the Laplacian of $h$ we have:
\begin{equation}\label{eq:est2}
|\Delta_\Sigma h|\le c_2(x,y)
\end{equation}
with $c_2(x,y)$ a smooth function that depends on $\tau$ and $H_0$.

Let $p$ be a point in $\H^2$ and $d$ denote the hyperbolic distance from $p$.
The function $d$ can be extended to the whole $\E(-1,\tau)$ by considering
$d\circ \pi$. We then have $\overline\nabla d=\wtilde{\nabla d}$ where $\wtilde
X$ denote the horizontal lift of $X$ (see Appendix~\ref{app:a}) and $\nabla$ is
the gradient operator on $\H^2$. We then have
\begin{equation}\label{eq:est3}
\|\nabla_\Sigma d^2\|=\|{\overline\nabla d^2}^\top\|\le d
\end{equation}
and because of formulas~\eqref{eq:levi1}, \eqref{eq:levi2} and \eqref{eq:levi3}
\begin{equation}\label{eq:est4}
|\Delta_\Sigma d^2|\le c_3(x,y,d)
\end{equation}
with $c_3$ a smooth function depending on $\tau$ and $H_0$


\subsection{A gradient estimate close to the boundary}

In this section, we give a first gradient estimate which is similar to Lemma 3.1
proved by Rosenberg, Schulze and Spruck in \cite{RoScSp}. This result will be
used to control the gradient of a section close to the boundary of the domain.
\begin{prop}\label{prop:graest1}
Let $\sigma$ be a section of $\E(-1,\tau)$ satisfying to \eqref{eq:cmc} which is
defined on a bounded domain $\Ome\subset\H^2$ and is $C^1$ up to the boundary.
Then there are two positive constants $M$ and $\alpha$ that depends only on
$\Ome$ and $\tau$ such that
$$
\sup_\Ome W\le \max(M,\sup_{\partial\Ome}W)\sup_\Ome e^{-\alpha \sigma}\sup_\Ome
e^{\alpha \sigma}.
$$
where $W=\sqrt{1+\|G\sigma\|^2}$. The quantity $e^{\pm\alpha \sigma}$ are computed
by identifying $\sigma$ with a function in the model.
\end{prop}

\begin{proof}
Let $\eta= e^{\alpha h}$ where $\alpha$ will be chosen later. From estimates
\eqref{eq:est1} and \eqref{eq:est2}, there is a positive constant $k$ that
depends only on $\Ome$ and $\tau$ such that
$$\Delta_\Sigma \eta=(\alpha^2 \|\nabla_\Sigma h\|^2+ \alpha\Delta_\Sigma
h)\eta\ge(\alpha^2(1-\nu^2-k\nu)-\alpha k)\eta$$
Thus 
$$L\frac\eta\nu=(L\frac1\nu)\eta+\frac1\nu\Delta_\Sigma\eta\ge
(-(1+2\tau^2) +\alpha^2(1-\nu^2-k\nu)-\alpha k)\frac\eta\nu$$

There is $\nu_0$ that depends only on $k$ such that, for $\nu\le \nu_0$:
$$-(1+2\tau^2) +\alpha^2(1-\nu^2-k\nu)-\alpha k\ge
-(1+2\tau^2)+\alpha^2/2-\alpha k$$
Thus, if $\alpha$ is chosen sufficiently large (depending only on $k$ and
$\tau$) we get $L\frac\eta\nu\ge 0$ where $\nu\le\nu_0$. By the maximum
principle, this implies that the maximum of $\frac\eta\nu$ is reached on the
boundary of $\Sigma$ or in $\{\nu\ge\nu_0\}$. This implies that
$$
\frac\eta\nu\le \max (\sup_{\partial\Sigma}\frac\eta\nu,
\frac1{\nu_0}\sup_\Sigma\eta)
$$
This is the expected estimate since $\nu=W^{-1}$.
\end{proof}


\subsection{A gradient estimate inside the domain}

We give now an other gradient estimate which is used to control a solution far
from the boundary. It is similar but less precise than the one given by Korevaar
in \cite{kor} or Spruck in \cite{Spr}

\begin{prop}\label{prop:graest2}
Let $\sigma$ be a section of $\E(-1,\tau)$ satisfying to \eqref{eq:cmc} which is
defined on a geodesic disk of $\H^2$ centered at $p=(x_p,y_p) \in D_{-1}$ and of
hyperbolic radius $R$. We assume that $\sigma$ viewed as a function in the model
is positive. Then there exists a positive constant $M$ that depends on
$x_p^2+y_p^2$, $R$, $\sigma(p)$ and $H_0$ such that 
$$
W(p)\le M
$$
\end{prop}

\begin{proof}
Let $d$ be the hyperbolic distance from $p$ in $\H^2$ and we extend it to the
whole $\E(-1,\tau)$. 

Let us define, on $\Sigma$, $\phi=(-\frac h{2h_0}+3/4-(\frac dR)^2)^+$ which
is less than $3/4$ and where $h_0=\sigma(p)$. If $P=(p,\sigma(p))$, we have
$\phi(P)=1/4$ and $\phi=0$ close to $\partial\Sigma$. Let us consider
$\eta=e^{K\phi}-1$ for a constant $K$ that will be chosen below. Let us define
$u=\frac\eta\nu$, we see that $\max u$ is positive and is reached inside the
support of $\phi$.

We have
\begin{align*}
Lu&=(L\frac 1\nu)\eta+(\frac 1\nu)\Delta_\Sigma \eta\\
&\ge -(1+2\tau^2)\frac\eta\nu+\frac 1\nu(K^2\|\nabla_\Sigma\phi\|^2+
K\Delta_\Sigma\phi)e^{K\phi}\\
&\ge (K^2\|\nabla_\Sigma\phi\|^2+
K\Delta_\Sigma\phi-(1+2\tau^2))\frac{e^{K\phi}}\nu +\frac{1+2\tau^2}\nu\\
&\ge (K^2\|\nabla_\Sigma\phi\|^2+ K\Delta_\Sigma\phi-(1+2\tau^2))
\frac{e^{K\phi}}\nu
\end{align*}
Let us see that $K$ can be chosen such that the first factor in the above
expression is positive.

Let us estimate the different terms. We have
$$
\nabla_\Sigma \phi=-\frac{\zeta^\top}{2h_0}-2\frac{\nabla_M d^2}{R^2}
$$
Thus, by \eqref{eq:zeta} and \eqref{eq:est3}, there is a constant $k_1$ that
depends on $p$ and $R$ such that $$
(\nabla_\Sigma \phi,\xi)\le-\frac{1}{2h_0}(1-k_1\nu)+k_1\nu
$$
So there exists $\nu_1>0$ that depends on $k_1$ and $h_0$ such that if
$\nu\le\nu_1$
$$
(\nabla_\Sigma \phi,\xi)\le-\frac{1}{2\sqrt{2}h_0}\textrm{\ \ and\ \ 
}\|\nabla_\Sigma \phi\|^2\ge \frac 1{8h_0^2}
$$

Besides, by \eqref{eq:est2} and \eqref{eq:est4}, there is a constant $k_2$ that
depends only on the domain $p$, $R$, $h_0$ such that $|\Delta_\Sigma \phi|\le
k_2$. This implies that for $\nu\le \nu_1$ we have
$$
K^2\|\nabla_\Sigma\phi\|^2+ K\Delta_\Sigma\phi-(1+2\tau^2)\ge
K^2\frac{1}{8h_0^2}-Kk_2-(1+2\tau^2)
$$
So there exists $K>0$ that depends only on $h_0$, $k_2$ and $\tau$ such that,
for $\nu\le\nu_1$,
$$
K^2\|\nabla_\Sigma\phi\|^2+ K\Delta_\Sigma\phi-(1+2\tau^2)>0
$$

Thus applying the maximum principle for the operator $L$, we get that the
maximum of $u$ is reached at a point $Q$ where $\nu\ge \nu_1$. Thus
$$
\frac{e^{K/4}-1}{\nu(P)}\le u(Q)\le\frac{e^{3K/4}-1}{\nu_1}
$$
So we get the expected estimate:
$$
\frac 1{\nu(P)}\le \frac 1{\nu_1}\frac{e^{3K/4}-1}{e^{K/4}-1}
$$
\end{proof}

\section{The half-space theorem}
\label{sec:half}
In this section, using Proposition~\ref{prop:uniq}, we prove the half space
theorem with respect to constant mean curvature $1/2$ entire graph in
$\E(-1,\tau)$. The theorem is the following.

\begin{thm}\label{th:half}
Let $\Sigma$ be a constant mean curvature $1/2$ entire graph in $\E(-1,\tau)$.
Let $\Sigma'$ be a properly immersed constant mean curvature $1/2$ surface in
$\E(-1,\tau)$ such that $\Sigma\cap\Sigma'=\emptyset$.
\begin{itemize}
\item if $\Sigma'$ is above $\Sigma$ then $\Sigma'=\Sigma+t$ for some $t>0$.
\item if $\Sigma'$ is below $\Sigma$ and is well oriented with respect to
$\Sigma$ then $\Sigma'=\Sigma-t$ for some $t>0$.
\end{itemize}
\end{thm}
Being well oriented means that the mean curvature vector of $\Sigma'$ points in
the connected component of $\E(-1,\tau)$ bounded by $\Sigma$ and $\Sigma'$.
Since $\Sigma'$ is only immersed, this condition has a meaning only for points
of $\Sigma'$ lying on the boundary of the connected component.

In the following, we write the proof only for the first case. We will just make
remarks to explain where the orientation hypothesis is used in the second case.
The approach is very similar to the one used by Rosenberg, Schulze and Spruck in
\cite{RoScSp} (see also Daniel, Meeks and Rosenberg \cite{DaMeRo}).


\subsection{Construction of the barriers}

We are working in the model. Let us consider an increasing sequence
$0<r_0<r_1<\cdots<r_n<\cdots$ with $r_n<1$ and $\lim r_n=1$ and let $A_n$ be the
annulus $\{r_0\le r\le r_n\}\subset D_{-1}$ and $A=\{r_0\le r<1\}$
($r^2=x^2+y^2$). We denote by $\sigma$ the function that defines the graph
$\Sigma$.

Using the implicit function theorem, there is $\delta>0$ such that there is a
smooth family $(u_{t,1})_{0\le t\le \delta}$ of smooth solutions $u_{t,1}$ of
\eqref{eq:cmc12} on $A_1$ with $u_{t,1}=\sigma+t$ on $\{r=r_0\}$ and
$u_{t,1}=\sigma$ on $\{r=r_1\}$.

In fact, the same construction can be made for every $n$ with the same $\delta$. 

\begin{lem}\label{lem:exist}
There is $\delta>0$ such that, for every $n\ge 1$ and $0\le t \le \delta$, there
exists a smooth solution $u_{t,n}$ of \eqref{eq:cmc12} on $A_n$ such
that $u_{t,n}=\sigma+t$ on $\{r=r_0\}$ and $u_{t,n}=\sigma$ on $\{r=r_n\}$. 

Moreover, for any $k>0$ the family of functions $(u_{t,n})_{0\le t\le \delta,n\ge
1}$ is uniformly bounded in the $C^{2,\alpha}$ norm on $A_k$.
\end{lem}

Because of the maximum principle, we have uniqueness of solutions to a Dirichlet
problem on compact domains; so the above solutions $u_{t,n}$ are unique. Thus
$u_{0,n}=\sigma$.
\begin{proof}
Let $\delta$ be given by the construction of $u_{t,1}$. Let us prove that this
constant works also for other $n$. Let us define $u_{0,n}=\sigma$ on $A_n$, the
existence of $u_{t,n}$ will come from the method of continuity. In order to
apply this method, we need some \textit{a priori} estimates for the solutions.

So let us consider a solution $u_{t,n}$ ($t\le\delta$). By the maximum
principle, we have $\sigma\le u_{t,n}\le \sigma+t\le \sigma+\delta$, so there is
a constant $c_1(n)>0$ such that $|u_{t,n}|\le c_1(n)$ on $A_n$. 

By the maximum principle, on $A_1$, we have $u_{\delta,1}+t-\delta\le
u_{t,n}\le \sigma+t$. But these three functions coincide on $\{r=r_0\}$ so the
gradient of $u_{t,n}$ on $\{r=r_0\}$ is bounded by a constant $c_2$ that does
not depend on $n$ and $t$.

From Appendix~\ref{app:b}, we know that there exists a smooth function $h$ on
$A_n$ such that $h=\sigma$ on $\{r=r_n\}$, $h\ge c_1(n)$ on $\{r=r_0\}$ and 
$$
\Div_{\H^2}\left(\frac{Gh}{\sqrt{1+\|Gh\|^2}}\right)\le 1
$$
so by the maximum principle, $\sigma\le u_{t,n}\le h$ on $A_n$. These functions
coincide on $\{r=r_n\}$ so the gradient of $u_{t,n}$ on $\{r=r_n\}$ is bounded
by a constant $c_3(n)$. Thus by Proposition~\ref{prop:graest1}, there is a
constant $c_4(n)$ such that $\|\nabla u_{t,n}\|\le c_4(n)$ on $A_n$ for
any $t\in[0,\delta]$. This implies that the equation solved by $u_{t,n}$ is
uniformly elliptic.

Then the DeGiorgi-Nash-Moser and Schauder estimates implies \textit{a priori}
bounds for higher derivatives of $u_{t,n}$ on $A_n$. The method of
continuity can then be applied to prove the existence of the solutions. We
notice that the estimates we just get depend on $n$ but in the lemma we want
estimates that are also independent of $n$.

For the $C^1$ estimate, let us consider $k\ge 1$, as above $\sigma\le u_{t,n}\le
\sigma+\delta$ so the family $(u_{t,n})_n$ is uniformly bounded on
$A_{k+1}$. Thus by Proposition~\ref{prop:graest2}, the gradient of
$u_{t,n}$ is uniformly bounded on $\{r=r_k\}$. Since $\|\nabla u_{t,n}\|\le c_2$
on $\{r=r_0\}$; Proposition~\ref{prop:graest1} tells that $\nabla u_{t,n}$ is
uniformly bounded on $A_k$. As above this gives uniform estimates for higher
derivatives on $A_k$.
\end{proof}

Since $(u_{\delta,n})$ and their derivatives are uniformly bounded, a diagonal
process gives a smooth solution $u_\delta$ of \eqref{eq:cmc12} on $A$
and a subsequence $(u_{\delta,n'})$ such that $u_{\delta,n'}\rightarrow
u_\delta$ where the convergence is smooth in all compact subsets of $A$.

By construction, we have $u_\delta=\sigma+\delta$ on $\partial A$ and $\sigma\le
u_\delta\le \sigma+\delta$. We remark that $\sigma+\delta$ is an other solution
of \eqref{eq:cmc12} with the same property; so by Proposition~\ref{prop:uniq} (see Section~\ref{sec:unik}),
$u_\delta=\sigma+\delta$.

\begin{remarq} For the case $\Sigma'$ below $\Sigma$, we need to construct
solutions $(u_{-\delta,n})$ which are below $\sigma$. The proof is similar, we
just use the $k$ barriers of Appendix~\ref{app:b} instead of the $h$ barriers.
\end{remarq}


\subsection{Proof of the half space theorem}

Let $\Sigma'$ be a properly immersed constant mean curvature $1/2$ surface in
$\E(-1,\tau)$ which is above $\Sigma$. By replacing $\Sigma$ by $\Sigma+t$ if
necessary, we can assume that $\Sigma'$ is not above $\Sigma+\eps$ for any
$\eps>0$.

If $\Sigma$ and $\Sigma'$ touches, the maximum principle implies that
$\Sigma'=\Sigma$; so we can assume that the two surfaces do not meet. 
Let $\delta'>0$ be such that $\Sigma'$ does not meet $\Sigma+t$ for
$t\in[0,\delta']$ over $\{r\le r_0\}$. Let $\delta\le \delta'$ such that
Lemma~\ref{lem:exist} is true and let us consider a subsequence of
$(u_{\delta,n})$ such that $u_{\delta,n'}\rightarrow \sigma+\delta$.

Let us denote by $\Sigma_{\delta,n}$ the graph of $u_{\delta,n}$. The surface
$\Sigma_{\delta,n}-\delta$ is below $\Sigma$ so below $\Sigma'$. Besides the
boundary of $\Sigma_{\delta,n}-t$ for $t\in[0,\delta]$ never meet $\Sigma'$,
thus, by the maximum principle, $\Sigma_{\delta,n}$ is below $\Sigma'$. Letting
$n$
tends to $+\infty$ along the chosen subsequence implies that $\Sigma+\delta=\lim
\Sigma_{\delta,n'}$ is below $\Sigma'$. This gives a contradiction and
Theorem~\ref{th:half} is proved.

\begin{remarq}
When $\Sigma'$ is below $\Sigma$, the orientation hypothesis is used in order to
apply the maximum principle between $\Sigma'$ and $\Sigma$ and between $\Sigma'$
and $\Sigma_{-\delta,n}+t$.
\end{remarq}

\section{A uniqueness exterior result}
\label{sec:unik}
In this section we prove a uniqueness result for \eqref{eq:cmc12} in an exterior
domain on $\H^2$.

We still consider the model for $\E(-1,\tau)$. Let us consider $r_0\in(0,1)$ and
$A=\{r\ge r_0\}\in D_{-1}=\H^2$. We have the following uniqueness result.
\begin{prop}\label{prop:uniq}
Let $u,v$ be two smooth solutions of \eqref{eq:cmc12} on $A$ such that $u=v$ on
$\partial A$ and $|u-v|$ is bounded on $A$ then $u=v$.
\end{prop}

The rest of this section is devoted to the proof of this result.

Let $u$ and $v$ be as in the proposition. If $u\neq v$ and exchanging $u$ and
$v$ if necessary, there is a $t_0>0$ such that $v+t_0\ge u$ and
$\inf_A(v-u+t_0)=0$. Let us denote $\tilde u=v+t_0$. We have $\tilde u = u+t_0$
on $\partial A$ and $\inf_A(\tilde u-u)=0$.


\subsection{Construction of a foliation}\label{sec:foliate}

Let $r_0<r_1<\cdots$ be an increasing sequence with $\lim r_n=1$. As above, we
denote $A_n=\{r_0\le r\le r_n\}$

Let $\delta<t_0$ be given by Lemma~\ref{lem:exist} with $\sigma=u$. Then we have
the associated family of solutions $u_{t,n}$ ($t\in[0,\delta]$) of
\eqref{eq:cmc12}.

We remark that if $t\le t'$, the maximum principle tells that $u_{t,n}\le
u_{t',n}\le u_{t,n}+t'-t\le \tilde u$. By
Lemma~\ref{lem:exist}, for any $k$, the family $(u_{t,n})_{0\le t\le \delta,
n\ge 1}$ is uniformly bounded in the $C^{2,\alpha}$ norm over $A_k$.

Thus by a diagonal process, there is a family of smooth solutions
$(u_{q\delta})_{q\in\Q\cap[0,1]}$ of \eqref{eq:cmc12} in $A$ and a
subsequence such that $u_{q\delta,n'}\rightarrow u_{q\delta}$ for any
$q\in\Q\cap[0,1]$ (the convergence is smooth in any compact subset of $A$).
If $q\le q'\in\Q\cap[0,1]$, we have $u_{q\delta}\le u_{q'\delta}\le
u_{q\delta}+(q'-q)\delta$. Besides the family $(u_{q\delta})_{q\in\Q\cap[0,1]}$
is uniformly bounded in the $C^{2,\alpha}$ norm over $A_k$. These
two properties imply that, for $t\in[0,\delta]$, the function 
$$
u_t=\sup_{\substack{q\in\Q\cap[0,1]\\ q\delta\le t}}u_{q\delta}=
\inf_{\substack{q\in\Q\cap[0,1]\\ q\delta\ge t}}u_{q\delta}
$$
is well defined and is a smooth solution of \eqref{eq:cmc12} on $A$ (we
notice that this definition coincide with the original one when
$t/\delta\in\Q$). By construction, if $t\le t'$, $u_t\le u_{t'}\le u_t+(t'-t)$
and $u_t=u+t=u_0+t$ on $\partial A$. Moreover, for any $k$, the family
$(u_t)_{t\in[0,\delta]}$ is uniformly bounded in the $C^{2,\alpha}$ norm
over $A_k$. Thus $\lim_{t\rightarrow t'}u_t=u_{t'}$ with smooth
convergence on any compact subset of $A$.  

We also have $u=u_0\le u_t\le u_\delta\le \tilde{u}$; this implies that
$\inf_A(u_t-u)=0$. 


\subsection{Derivation of the foliation}

In the preceding subsection, we have constructed the family $(u_t)$ and proved
that it depends continuously in the parameter $t$. In this subsection, we study
the differentiability with respect to $t$. 

Since all the functions $u_t$ satisfies to \eqref{eq:cmc12}, a derivative
$v=\frac{du_t}{dt}|_{t=\bar t}$ is \textit{a priori} a solution of
\begin{equation}\label{eq:der}
\Div_{\H^2}\left(\frac{\nabla v-\chi_{u_{\bar{t}}}(\nabla v,\chi_{u_{\bar t}})}
{\sqrt{1+\|Gu_{\bar t}\|^2}}\right)=0
\end{equation}
where $\dis \chi_u=\frac{Gu}{\sqrt{1+\|Gu\|^2}}$ (see Appendix~\ref{app:a}). We
will prove that this is in fact the case.

\begin{lem}\label{lem:deriv}
Let $\bar t\in[0,\delta]$ and $(t_k)$ a sequence converging to $\bar t$ then
there is a solution $v$ of \eqref{eq:der} and a subsequence such that 
$$
\frac{u_{t_{k'}}-u_{\bar t}}{t_{k'}-\bar t}\longrightarrow v
$$
with a smooth convergence in any compact subset of $A$.
\end{lem}

The idea of the proof comes from the work of Solomon about the regularity of
minimal foliations \cite{Sol}.

\begin{proof}
Let $v_k$ be equal to $u_{t_k}-u_{\bar t}$. We then have $Gu_{t_k}=Gu_{\bar
t}+\nabla v_k$. Thus $v_k$ is a solution of the following equation
\begin{align}
0&=\Div_{\H^2}\left(\frac{\nabla v_k+Gu_{\bar t}} {\sqrt{1+\|\nabla v_k+
Gu_{\bar t}\|^2}}-\frac{Gu_{\bar t}} {\sqrt{1+\|Gu_{\bar t}\|^2}}\right)\notag\\
&=\Div_{\H^2}\left(\int_0^1 \frac{\nabla v_k-\chi_{w_{k,s}}(\nabla
v_k,\chi_{w_{k,s}})} {\sqrt{1+\|Gw_{k,s}\|^2}} \dd s \right)
\end{align}
with $w_{k,s}=u_{\bar t}+sv_k$. Let us denote by $P_{k,s}$ the linear operator 
$$
P_{k,s}(X)=\frac {X-\chi_{w_{k,s}}(X,\chi_{w_{k,s}})} {\sqrt{1+\|Gw_{k,s}\|^2}}.
$$
Actually, $P_{k,s}$ is a smooth section of the endomorphisms on $T\H^2$ defined
on $A$. 

Let $n\in\N^*$, we know that the functions $u_{t_k}$ and $u_{\bar t}$ are
uniformly bounded in $C^{2,\alpha}$ norm over $A_n$. Thus the family
$(w_{k,s})_{k\in\N,s\in[0,1]}$ is uniformly bounded in $C^{2,\alpha}$ norm over
$A_n$. This implies first that the family $(P_{k,s})$ is uniformly bounded in
$C^{1,\alpha}$ norm over $A_n$ and that there exist a constant $c_n>0$ such that,
for any $X\in T\H^2$,
$$
(P_{k,s}X,X)=\frac{\|X\|-(X,\chi_{w_{k,s}})^2}{\sqrt{1+\|Gw_{k,s}\|^2}}\ge
\frac{(1-\|\chi_{w_{k,s}}\|^2)\|X\|^2}{\sqrt{1+\|Gw_{k,s}\|^2}}\ge c_n\|X\|^2
$$
Thus we have proved that the operators
$$
L_k(w)=\Div_{\H^2}\int_0^1P_{k,s}(\nabla w)\dd s
$$
are uniformly elliptic with uniformly bounded $C^{0,\alpha}$ coefficients in
$A_n$. So, by Schauder estimates, there is a constant $c_n$ such that
$$
\|v_k\|_{C^{2,\alpha}(A_n)}\le c_n(\|v_k\|_{C^0(A_{n+1})}+
\|v_k\|_{C^{2,\alpha}(\{r=r_0\})}) \le 2c_n|t_k-\bar t|
$$
The last inequality comes from the fact that $|v_k|=|u_{t_k}-u_{\bar t}|\le
|t_k-\bar t|$ and $v_k=t_k-\bar t$ along $\{r=r_0\}$. Thus $v_k/(t_k-\bar t)$ is
uniformly bounded in $C^{2,\alpha}$ norm over $A_n$. So, by a diagonal process,
there is a subsequence $v_{k'}/(t_{k'}-\bar t)$ that converges smoothly to
a function $v$ on any compact subset of $A$.

Besides, since $v_k\rightarrow 0$ smoothly, the operator $P_{k,s}$ converges
uniformly to the operator $P$ defined by:
$$
P(X)=\frac {X-\chi_{u_{\bar t}}(X,\chi_{u_{\bar t}})} {\sqrt{1+\|Gu_{\bar
t}\|^2}}.
$$
So $v$ is a solution of \eqref{eq:der}.
\end{proof}

Since $u_t\le u_{t'}\le u_{t}+(t'-t)$ if $t\le t'$, the function $v$ satisfies
$0\le v\le 1$. Moreover, on $\partial A$, we have $v=1$. The constant function
$1$ is an obvious solution to \eqref{eq:der} that satisfies to the same
properties; in the following we will indeed prove that $1$ is the only such
solution.


\subsection{The associate surface in $\E(0,\tau')$}

Let $\bar t$ be in $[0,\delta]$ and $\Sigma$ be the graph of $u_{\bar t}$. We
will consider several metrics on the annulus $\Sigma$, so in the computations, we
will always explain to which metric the computation is made. The function $v$
constructed by Lemma~\ref{lem:deriv} can be viewed as a function on
$\Sigma$. From Appendix~\ref{app:a}, if $\nu$ is the angle function $(N,\xi)$,
$v\nu$ is a Jacobi function on $(\Sigma,g)$
where $g$ is the induced metric from $\E(-1,\tau)$. In this section, we use the
associated family of constant mean curvature surfaces introduced by Daniel
\cite{Dan2} and the work of Daniel and Hauswirth \cite{DaHa} to prove the
following result.

\begin{lem}\label{lem:flat1}
There exist a flat metric $g_0$ on $\Sigma$ and a vector field $G$ on
$\Sigma$ such that any function $w$ such that $w\nu$ is a Jacobi function of
$(\Sigma,g)$ is a solution to 
\begin{equation}\label{eq:jacobquotient}
\Div_0\left(\frac{\nabla^0 w-\chi(\nabla^0 w,\chi)}
{\sqrt{1+\|G\|_0^2}}\right)=0
\end{equation}
with $\chi=\frac G{\sqrt{1+\|G\|_0^2}}$ (the sub or superscript $0$ means that
the computation are made with respect to the $g_0$ metric).
\end{lem}

\begin{proof}
Let $\widetilde\Sigma$ be the universal cover of $\Sigma$. The immersion $X$ of
$\widetilde\Sigma$ in $\E(-1,\tau)$ is encoded in the induced metric $g$, the
shape operator $S$, the function $\nu$ and the orthogonal projection $T$ of the
vector field $\xi$. As in Section~\ref{sec:entiregraph}, let $\theta$ be such
that $\tau+\frac i2= e^{i\theta} \sqrt{\frac 14+\tau^2}$. Then, the data
\begin{align*}
g'&=g\\
S'&=e^{\theta J}(S-\frac 12I)\\
\nu'&=\nu\\
T'&=e^{\theta J}T
\end{align*}
encode a minimal immersion $X'$ of $\widetilde\Sigma$ in $\E(0,\sqrt{\frac
14+\tau^2})$ ($J$ denotes the rotation by $\pi/2$ in the tangent space to
$\widetilde\Sigma$) (see \cite{Dan2}).

The Jacobi operator of $X'(\widetilde\Sigma)$ is:
\begin{align*}
\Delta_{g'}+(-2\tau'^2+4\tau'^2\nu'^2+\|S'\|^2)&=\Delta_g+ (-2(\frac 14+\tau^2)+
4(\frac 14+\tau^2)\nu^2+\|S\|^2-\frac 12)\\
&=\Delta_g+(-\frac 12-2\tau^2)+(1+4\tau^2)\nu^2+\|S\|^2)
\end{align*}
where the computation are made with respect to the metric $g=g'$.
Thus the Jacobi operator of $X'(\widetilde\Sigma)$ is the same that the one of
$X(\widetilde\Sigma)$. If $w$ is viewed as a function on $\widetilde\Sigma$, the
function $w\nu$ is a Jacobi function on $X(\widetilde\Sigma)$ so it is a Jacobi
function on $X'(\widetilde\Sigma)$.

Let $\pi$ be the submersion from $\E(0,\tau')$ to $\R^2$. Since $\nu'=\nu>0$,
the map $\pi\circ X'$ is a local diffeomorphism from $\widetilde\Sigma$ to $\R^2$,
so we can lift to $\widetilde\Sigma$ the flat metric of $\R^2$. Let $g_0$ denote
this flat metric on $\widetilde\Sigma$. Besides locally, we can describe
$X'(\widetilde\Sigma)$ has the graph of a section $s$. By Appendix~\ref{app:a},
since $w\nu'$ is a Jacobi function on $X'(\widetilde\Sigma)$, the function $w$,
viewed as a function on $\R^2$, is a solution of 
$$
\Div_{\R^2}\left(\frac{\nabla^{\R^2} w-\chi_s(\nabla^{\R^2} w,\chi_s)}
{\sqrt{1+\|Gs\|_{\R^2}^2}}\right)=0
$$
The computation are made with respect to the Euclidean metric of $\R^2$.
But the vector field $G=(\pi\circ X')^*(Gs)$ is globally well defined on
$\widetilde\Sigma$ and then $v$ is a solution on $\widetilde\Sigma$ of:
$$
\Div_0\left(\frac{\nabla^0 w-\chi(\nabla^0 w,\chi)}
{\sqrt{1+\|G\|_0^2}}\right)=0
$$
where $\chi=\frac G{\sqrt{1+\|G\|_0^2}}$. Let us now see that this description
passes to the quotient surface $\Sigma$. Let us consider $\gamma$ a generator of
$\pi_1(\Sigma)$. The element $\gamma$ acts on $\widetilde\Sigma$ as a
diffeomorphism without any fixed point. Moreover, from the uniqueness part of
Theorem 4.3 in \cite{Dan2}, there exists $M$ an isometry of $\E(0,\tau')$ such
that for any $p\in\widetilde\Sigma$ we have $X'(\gamma\cdot p)=M(X'(p))$.
Besides there is an isometry $m$ of $\R^2$ such that $\pi(M(q))=m(\pi(q))$ for
any $q\in\E(0,\tau')$. This implies that $\pi(X'(\gamma\cdot p))=m(\pi(X'(p)))$
so $\gamma^*g_0=g_0$ and $\gamma^* G=G$. This implies that the metric and the
vector field pass to the quotient surface $\Sigma$ and $w$ on $\Sigma$ is a
solution of \eqref{eq:jacobquotient}.
\end{proof}

In fact, the metric $g_0$ on the surface $\Sigma$ which has a compact boundary
satisfies the following property.

\begin{lem}\label{lem:flat2}
The flat metric $g_0$ is complete.
\end{lem}

\begin{proof}
It suffices to prove that every curve in $\Sigma$ starting from a point
in the boundary of $\Sigma$ which is proper in $\Sigma$ has infinite length with
respect to $g_0$. On $\Sigma$, we lift the function $r$ which is the radial
coordinate on $D_{-1}$; $r:\Sigma\rightarrow[0,1)$ is proper. $r$ can also be
lifted to $\widetilde\Sigma$.

Actually, it suffices to prove that every curve in $\widetilde\Sigma$ starting
from a point in the boundary of $\widetilde\Sigma$ such that $r\circ
\gamma\rightarrow 1$ has infinite length with respect to $g_0$.

So let $\gamma$ be such a curve in $\widetilde\Sigma$. We have the following
estimate:
$$
\ell(\gamma,g_0)=\ell(\pi(X'(\gamma)),\R^2)\ge \int_{X'(\gamma)}\nu'\dd
\ell_{g'}= \int_{X(\gamma)}\nu\dd \ell_{g} \ge \int_{\pi(X(\gamma))}\nu \dd
\ell_{\H^2}
$$
In $\{r_0\le r\le r_1\}$, we know that there is $\nu_0>0$ such that
$\nu\ge \nu_0$. So for a curve $\gamma$ with $r\circ \gamma$ proper, we have:
$$
\ell(\gamma,g_0)\ge 2\nu_0(\argth r_1-\argth r_0)=\ell_0>0
$$
This estimate implies that we can only consider points in $\widetilde\Sigma$
which are at a distance larger than $\ell_0$ from $\partial\widetilde\Sigma$ in
the $g_0$ metric. Using this, we can apply the work of Daniel and
Hauswirth~\cite{DaHa} to prove that, if $(\widetilde\Sigma,g_0)$ does not
satisfied to the expected property, then $X'(\widetilde\Sigma)$ has a subset
which is a graph over a strip $S\in\R^2$ isometric to $(0,\eps)\times\R$.
Moreover this graph goes to $+\infty$ on one of its boundary component. But the
existence of such a minimal graph is impossible by Theorem 6.3 in \cite{DaHa}.
\end{proof}


\subsection{A uniqueness result for Jacobi function}

In this subsection, we prove a uniqueness result for solution of a certain
partial derivative equation over a flat surface.

\begin{lem}\label{lem:uniqueness}
Let $(S,\dd s^2)$ be a complete flat surface with a compact boundary. Let $G$ be
a smooth vector field on $S$. We consider $w$ and $w'$ two smooth functions on
$S$ with the same boundary value which solve the following partial derivatives
equation:
$$
\Div\left(\frac{\nabla u-\chi(\nabla u,\chi)}
{\sqrt{1+\|G\|^2}}\right)=0
$$
with $\chi=\frac G{\sqrt{1+\|G\|^2}}$. If $w-w'$ is bounded then $w=w'$.
\end{lem}

\begin{proof}
The first step of the proof is to estimate the growth of the surface $S$. Let
$d$ be the distance function from $\partial S$. $d$ is a Lipschitz function. If
$n$ is the inward unit normal vector to $\partial S$, the set $\{d=d_0\}$ is
included in the image of the map $\phi_{d_0}:\partial S\rightarrow S; p\mapsto
\exp_p(d_0N(p))$ (actually this map is well defined only on a subset of
$\partial S$). With $\kappa$ the geodesic curvature of $\partial S$, we obtain:
$$
\ell(\{d=d_0\}\le \ell(\im \phi_{d_0})\le\int_{\partial S}|1+\kappa d_0|\dd\ell
\le (1+\bar\kappa d_0)\ell(\partial S)
$$
where $\bar\kappa=\max |\kappa|$. Thus $\ell(\{d=d_0\}$ has at most a linear
growth.

Let $P$ be the linear map defined on $TS$ by:
$$
P(X)=\frac{X-\chi(X,\chi)}{\sqrt{1+\|G\|^2}}
$$
We have 
\begin{align*}
\|P(X)\|^2&=\frac{\|X\|^2-(X,\chi)^2}{1+\|G\|^2}-
\frac{(1-\|\chi\|^2)(X,\chi)^2} {1+\|G\|^2}\\
&\le \frac{\|X\|^2-(X,\chi)^2}
{1+\|G\|^2}= \frac{(X,P(X))}{\sqrt{1+\|G\|^2}}\\
&\le(X,P(X))
\end{align*}

The function $d$ is Lipschitz continuous and the set $\{0\le d\le d_0\}$ is
compact with rectifiable boundary. So let us define 
$$
I(d_0)=\int_{\partial\{0\le d\le d_0\}}(w-w')(P(\nabla w-\nabla
w'),\eta)\dd \boH^1
$$
and $\dis \mu(d_0)= \int_{\{d=d_0\}}\|P(\nabla w-\nabla w')\|\dd \boH^1$. Since
$w-w'=0$ on $\{d=0\}$, $|w-w'|\le M$ and $\partial\{0\le d\le d_0\}\subset
\{d=0\}\cup\{d=d_0\}$, we have $|I(d_0)|\le M \mu(d_0)$. Using Stokes and coarea
formulas and $\|\nabla d\|=1$ a.e., we also have for $d_0>d_1$:
\begin{align*}
I(d_0)&=\int_{\{d\le d_0\}}(\nabla w-\nabla w',P(\nabla w-\nabla w'))\\
&\ge I(d_1)+\int_{\{d_1\le d\le d_0\}}\|P(\nabla w-\nabla w')\|^2\\
&\ge I(d_1)+\int_{d_1}^{d_0}\int_{\{d=s\}}\|P(\nabla w-\nabla
w')\|^2\dd\boH^1\dd s
\end{align*}
Since $\ell(\{d=s\})$ is at most linear, there is a $c>0$ such that for $s\ge
d_1$:
$$
\mu(s)^2\le \ell(\{d=s\})\int_{\{d=s\}}\|P(\nabla w-\nabla
w')\|^2\dd\boH^1 \le c s \int_{\{d=s\}}\|P(\nabla w-\nabla
w')\|^2\dd\boH^1
$$
So we have:
$$
\frac{I(d_1)}{M}+\int_{d_1}^{d_0}\frac{\mu^2(s)}{c M s}\dd s \le \mu(d_0)
$$
We notice that $\mu$ is \textit{a priori} just a locally bounded measurable
function, but we can integrate the above differential inequality. Indeed let
$g(d_0)$ be the left-hand side of the above inequality. Since $\mu$ is locally
bounded, $g$ is locally Lipschitz. Moreover
$$
g(d_0)=\frac{I(d_1)}{M}+\int_{d_1}^{d_0}\frac{\mu^2(s)}{c M s}\dd s\ge
\frac{I(d_1)}{M}+\int_{d_1}^{d_0}\frac{g^2(s)}{c M s}\dd s
$$
Now if $f(d_0)$ is the right-hand side of the above inequality, $f$ is $C^1$,
satisfies the same differential inequality and we have :
$$
f'(s)=\frac{g^2(s)}{c M s}\ge \frac{f^2(s)}{c M s}
$$
So 
$$
f(d_0)\ge \frac {I(d_1)}{1-\frac {I(d_1)}{M^2c}\log(\frac{d_0}{d_1})}
$$
Now since $f(d_0)\le g(d_0) \le \mu(d_0)$ 
$$
\mu(d_0)\ge \frac {I(d_1)}{1-\frac {I(d_1)}{M^2c}\log(\frac{d_0}{d_1})}
$$
But this lower bound blows up at $d_0=d_1\exp(\frac{M^2c}{I(d_1)})$, so we
get a contradiction and this finishes the proof of Lemma~\ref{lem:uniqueness}.
\end{proof}


\subsection{End of the proof}

Let us now use all the preceding lemmas to finish the proof of
Proposition~\ref{prop:uniq}.

So let $u$ and $\tilde u$ be the function introduced at the beginning of the
section. Let $(u_t)_{t\in[0,\delta]}$ be the family of function constructed in
Section~\ref{sec:foliate}. 

Let $\bar t\in[0,\delta]$ and $(t_k)$ be a sequence converging to $\bar t$. From
Lemma~\ref{lem:deriv}, there are a solution $v$ of \eqref{eq:der} and a
subsequence such that
$$
\frac{u_{t_{k'}}-u_{\bar t}}{t_{k'}-\bar t}\longrightarrow v
$$
By construction $0\le v\le 1$ and $v=1$ on $\partial A$. $w\equiv1$ is an other
solution of \eqref{eq:der} with the same properties. Let us see $v$ and $w$ as
functions on the graph $\Sigma$ of $u_{\bar t}$. From Lemmas~\ref{lem:flat1} and
\ref{lem:flat2}, there is a complete flat metric $g_0$ on $\Sigma$ and a vector
field $G$ such that $v$ and $w$ are solutions of \eqref{eq:jacobquotient}. Thus,
by Lemma~\ref{lem:uniqueness}, $v=w\equiv 1$. The uniqueness of the possible
limit implies that
$$
\frac{u_{t_k}-u_{\bar t}}{t_k-\bar t}\longrightarrow 1
$$
for every sequence $(t_k)$. So $u_t$ is differentiable with respect to $t$ and
$\der{u_t}{t}=1$. This implies that $\tilde u\ge u_\delta=u_0+\delta=u+\delta$. So
we have a contradiction with $\inf_A (\tilde u-u)=0$ and it finishes the proof
of Proposition~\ref{prop:uniq}

\appendix
\section{Some computations in Killing Riemannian submersions}
\label{app:a}

In this appendix we recall some definitions about Killing Riemannian submersions
and make some computations concerning graphs in such an ambient space (see
\cite{Man,EsOl}).

Let $(M^{n+1},\bar g)$ and $(B^n,(\cdot,\cdot))$ be two complete Riemannian
manifolds. Let
$\pi:M\rightarrow B$ be a submersion. The tangent space $T_pM$ at $p$ then
splits in $\ker d\pi\oplus (\ker d\pi)^\perp$ where $\ker d\pi$ is the
$1$-dimensional space of \textit{vertical} vectors and $(\ker d\pi)^\perp$ is
the space of \textit{horizontal} vectors. The submersion $\pi$ is called
\textit{Riemannian} if $d\pi$ is an isometry from $(\ker d\pi)^\perp$ to
$T_{\pi(p)}B$.

\begin{defn}
A Riemannian submersion $\pi:M\rightarrow B$ is a Killing submersion if it
admits a complete vertical unit Killing vector field.
\end{defn}

If $\pi:M\rightarrow B$ is a Killing Riemannian submersion, we denote by
$\xi$ this unit Killing vector field. Besides, if $X$ is a vector field in $B$,
we denote by $\widetilde X$ its horizontal lift by $\pi$.

Using this notation there exists a $2$-form $\omega$ on $B$ such that for any
vector fields $X,Y$ on $B$ we have:
$$
[\wtilde X,\wtilde Y]=\wtilde{[X,Y]}+\omega(X,Y)\xi
$$
We notice that $[\wtilde X,\xi]=0$. When $X$ is a tangent vector to $B$, we
denote $X^\omega$ the vector such that $(X^\omega,Y)=\omega(X,Y)$ for any $Y$.
With this notation the Levi-Civita connection of $M$ and $B$ are relied by 
\begin{align}\label{eq:levi1}
&\overline\nabla_{\wtilde X}\wtilde Y=\wtilde{\nabla_XY}+\frac12\omega(X,Y)\xi\\
\label{eq:levi2}
&\overline\nabla_{\wtilde X}\xi=-\frac12\wtilde{X^\omega}\\
\label{eq:levi3}
&\overline\nabla_\xi\xi=0
\end{align}

In a Killing Riemannian submersion, we are interested by surfaces that are the
image of sections $\sigma$. These surfaces are called vertical graphs in $M$ and
the image of $\sigma$ is also called the graph of $\sigma$. If
$\sigma$ is a section defined over $\Ome\subset B$, we define on $\Ome$ a vector
field $G\sigma$ by the following property:
$$
(G\sigma,X)=\bar g(\dd \sigma(X),\xi)
$$
for any $X$ tangent to $B$. This vector field $G\sigma$ plays the role of the
gradient of a function. 

First, the upward pointing unit normal to the graph of $\sigma$ is given by the
following expression:
$$
N=\frac{-\wtilde{G\sigma}+\xi}{\sqrt{1+\|G\sigma\|^2}}
$$
In the following, we denote $\sqrt{1+\|G\sigma\|^2}$ by $W$. In fact the
expression of $N$ is defined in the whole $\pi^{-1}(\Ome)$ and the mean
curvature of the graph of $\sigma$ is given by
$$
nH=-\Div_{\bar g}\left(\frac{-\wtilde{G\sigma}+\xi}{W}\right)
$$ 
So if $(E_i)$ is an orthonormal frame of $T_pB$ we have
\begin{align*}
nH&=-\sum_i\bar g(\overline\nabla_{\wtilde E_i} \frac{-\wtilde{G\sigma}+\xi}
{W}, \widetilde E_i)-\bar g(\overline\nabla_\xi \frac{-\wtilde{G\sigma}+\xi}
{W}, \xi)\\
&=\sum_i\bar g(\wtilde{\nabla_{E_i}\frac{G\sigma}W},\wtilde{E_i})+
\frac1{2W}\bar g(\wtilde{E_i^\omega},\wtilde{E_i})\\
&=\Div_g\left(\frac{G\sigma}W\right)+\sum_i\frac1{2W}\omega(E_i,E_i)\\
&=\Div_g\left(\frac{G\sigma}W\right)
\end{align*}
So the mean curvature is given by 
\begin{equation}\label{eq:appcmc}
nH=\Div\left(\frac{G\sigma}{\sqrt{1+\|G\sigma\|^2}}\right)
\end{equation}

Let us denote by $\Sigma$ the graph of $\sigma$. The map $\sigma:
\Ome\rightarrow \Sigma$ is a chart, so let us make some computation using this
system of coordinates. We have:
$$
\dd \sigma(X)=\widetilde{X}+(G\sigma,X)\xi
$$
so the induced metric is $g(X,Y)=(X,Y)+(G\sigma,X)(G\sigma,Y)$. If $u$ is a
function on $\Ome$, we get
$$
(\nabla u,X)=\dd u(X)=g(\nabla_g u,X)=(\nabla_g u,X)+(G\sigma,\nabla_g u)
(G\sigma,X)
$$
Thus $\nabla u=\nabla_g u+(G\sigma,\nabla_g u)G\sigma$; this implies that
$$
\nabla_g u=\nabla_u-(\chi_\sigma,\nabla u)\chi_\sigma
$$
where $\chi_\sigma=\frac{G\sigma}{\sqrt{1+\|G\sigma\|^2}}$. As a consequence, we
have:
\begin{align*}
g(\nabla_g v,\nabla_g w)&=(\nabla v-(\chi_\sigma,\nabla v)\chi_\sigma, \nabla
w-(\chi_\sigma,\nabla w)\chi_\sigma)\\
&\quad\quad\quad+(G\sigma, \nabla v-(\chi_\sigma, \nabla
v)\chi_\sigma)(G\sigma, \nabla w-(\chi_\sigma, \nabla w)\chi_\sigma)\\
&=(\nabla v,\nabla w)-(2-\|\chi_\sigma\|^2)(\chi_\sigma, \nabla v)(\chi_\sigma,
\nabla w)\\
&\quad\quad\quad+W^2(1-\|\chi_\sigma\|^2)^2(\chi_\sigma, \nabla v)(\chi_\sigma,
\nabla w)\\
&=(\nabla v,\nabla w)-(\chi_\sigma, \nabla v)(\chi_\sigma, \nabla w)
\end{align*}

Besides the divergence operator for the metric $g$ have the following
expression:
$$
\Div_g X=\frac1{\sqrt{1+\|G\sigma\|^2}}\Div\big(\sqrt{1+\|G\sigma\|^2} X\big)
$$

If $\Sigma$ has constant mean curvature, the function $\nu=(N,\xi)$ is a Jacobi
function so $0=\Delta_\Sigma \nu+(Ric_{\bar g}(N,N)+\|S\|^2)\nu$ with $S$ the
shape operator of $\Sigma$. Let $v$ be a function on $\Sigma$, we then have
\begin{align*}
\Delta_\Sigma (\nu v) +(Ric_{\bar g}(N,N)+\|A\|^2) (\nu v)&= (\Delta_\Sigma
\nu)v+2(\nabla_\Sigma u,\nabla_\Sigma v)+\nu(\Delta_\Sigma v)+ (Ric_{\bar
g}(N,N)+\|A\|^2) (\nu v)\\
&=\nu\big(\Delta_\Sigma v-2\nu(\nabla_\Sigma\frac{1}{\nu},\nabla_\Sigma v)\big)
\end{align*}
Thus $\nu v$ is a Jacobi function if and only if $\Delta_\Sigma
v-2\nu(\nabla_\Sigma\frac{1}{\nu},\nabla_\Sigma v)=0$. Looking at $v$ as a
function defined on $\Ome$, since $\nu=W^{-1}$, it gives:
\begin{align*}
0&=\Delta_g v-2W^{-1}g(\nabla_g W,\nabla_g v)\\
&=W^{-1}\big(\Div(W(\nabla v-(\chi_\sigma,\nabla v)\chi_\sigma))-2((\nabla
W,\nabla v)-(\chi_\sigma, \nabla W) (\chi_\sigma, \nabla v))\big)\\
&=W^{-1}\big(\Div(W^2\frac{\nabla v-(\chi_\sigma,\nabla
v)\chi_\sigma}{W})-(\nabla W^2, \frac{\nabla v-(\chi_\sigma,\nabla
v)\chi_\sigma}{W})\big)\\
&=W^{-1}\big(\Div(\frac{\nabla v-(\chi_\sigma,\nabla
v)\chi_\sigma}{W})\big)
\end{align*}

So $v\nu$ is a Jacobi function  if and only if
\begin{equation}
0=\Div(\frac{\nabla v-(\chi_\sigma,\nabla v)\chi_\sigma}{W})
\end{equation}


\section{Some barriers}
\label{app:b}

In this appendix, we construct some barriers from above and below on the exterior
boundary component of an annulus for \eqref{eq:cmc12}. 

We use the model for $\E(-1,\tau)$ but we consider hyperbolic polar coordinates
on $D_{-1}$, so $x=\tanh(\rho/2)\cos\theta$ and $y=\tanh(\rho/2)\sin\theta$.

Let $0<\rho_1<\rho_0$ be two radii and $f$ be a smooth function on
$\{\rho=\rho_0\}$. For $M$ a constant, we construct a smooth function
$h$ on $\{\rho_1\le \rho \le \rho_0\}$ such that 
\begin{itemize}
\item $h=f$ on $\{\rho=\rho_0\}$, $h\ge M$ in $\{\rho=\rho_1\}$ and 
\item $\dis \Div_{\H^2}\left(\frac{Gh}{\sqrt{1+\|Gh\|^2}}\right)\le 1$
\end{itemize}

We see $f$ as a function of $\theta$ and we define on $\{\rho_1\le \rho \le
\rho_0\}$ the function $h(\rho,\theta)=f(\theta)+\alpha(\rho_0-\rho)$. Let us
prove that for $\alpha$ sufficiently large $h$ satisfies the expected
properties. We have $h(\rho_0,\theta)=f(\theta)$ and
$h(\rho_1,\theta)=f(\theta)+\alpha(\rho_0-\rho_1)\ge M$ if $\alpha\ge (M-\min
f)/(\rho_0-\rho_1)$. In the polar coordinates, we have
$$
Gh=-\alpha\partial_\rho+ \frac1{\sinh\rho}\big(\frac{f'(\theta)}{\sinh\rho}-
2\tau\tanh\frac{\rho}2\big)\partial_\theta
$$
Thus
\begin{align*}
\Div_{\H^2}\left(\frac{Gh}{\sqrt{1+\|Gh\|^2}}\right)&=
-\frac{\cosh\rho}{\sinh\rho}\frac{\alpha}{W}+\alpha
\frac{\partial_\rho(\frac{f'(\theta)}{\sinh\rho}-
2\tau\tanh\frac{\rho}2\big)^2}{2W^3}+\frac1{\sinh\rho}
\frac{\partial_\theta(\frac{f'(\theta)}{\sinh\rho})}{W}\\
&\quad\quad -\frac1{\sinh\rho}
\frac{(\frac{f'(\theta)}{\sinh\rho}-
2\tau\tanh\frac{\rho}2\big)\partial_\theta(\frac{f'(\theta)}{\sinh\rho}-
2\tau\tanh\frac{\rho}2\big)^2}{2W^3}
\end{align*}
where
$$
W=\sqrt{1+\alpha^2+\big(\frac{f'(\theta)}{\sinh\rho}-
2\tau\tanh\frac{\rho}2\big)^2}
$$
With $\alpha>0$, there exists a positive constant $m$ that only depends on $f$,
$\rho_0$ and $\rho_1$ such that:
$$
\Div_{\H^2}\left(\frac{Gh}{\sqrt{1+\|Gh\|^2}}\right)\le
\frac{m}{\alpha^2}+\frac{m}\alpha+\frac{m}{\alpha^3}
$$
So when $\alpha$ is sufficently large, the mean curvature of the graph of $h$
satisfies the expected estimate.

Now let us define $k(\rho,\theta)=f(\theta)+\alpha(\rho_0-\rho)$ with $\alpha<0$
By choosing $\alpha$ small, we can ensure that $k(\rho_1,\theta)\le M$.
Moreover, because of the above computation, there is a $m>0$ such that
$$
\Div_{\H^2}\left(\frac{Gk}{\sqrt{1+\|Gk\|^2}}\right)\ge
\frac{\cosh\rho_0}{\sinh\rho_0}\frac{-\alpha}{\sqrt{\alpha^2+m}}-
\frac{m}{\alpha^2}-\frac{m}{|\alpha|}-\frac{m}{|\alpha|^3}
$$
So choosing $\alpha$ sufficiently close from $-\infty$, $k$ satisfies to
\begin{itemize}
\item $k=f$ on $\{\rho=\rho_0\}$, $k\le M$ in $\{\rho=\rho_1\}$ and 
\item $\dis \Div_{\H^2}\left(\frac{Gk}{\sqrt{1+\|Gk\|^2}}\right)\ge 1$
\end{itemize}

\bibliographystyle{plain}
\bibliography{../reference.bib}

\noindent \textsc{Laurent Mazet}\\
Universit\'e Paris-Est \\
Laboratoire d'Analyse et Math\'ematiques Appliqu\'ees, CNRS UMR8050\\
UFR des Sciences et Technologie\\ 
61 avenue du G\'en\'eral de Gaulle\\
94010 Cr\'eteil cedex, France

\noindent \texttt{laurent.mazet@math.cnrs.fr}

\end{document}